\documentclass[12pt]{article}

\usepackage{amsmath,latexsym,amssymb,amstext,amsthm}

\usepackage[dvips]{graphicx}
\usepackage[latin1]{inputenc}

\begin{document}

\title{How to model marine reserves ? }
\author{Patrice Loisel \footnote{e-mail : patrice.loisel@supagro.inra.fr,
phone : 33 (0) 4 99 61 29 04, fax : 33 (0) 4 67 52 14 27, address : INRA-UMR-ASB 2 place Viala 34060 MONTPELLIER FRANCE} Pierre Cartigny \footnote{e-mail : pierre.cartigny@supagro.inra.fr}}

\date{}

\maketitle

\baselineskip=16pt

%Article Type: Research Paper \\

\vspace{10mm}

Abstract: The safeguarding of resources is one of the principal subjects of halieutics studies. Among the solutions proposed to avert the disappearance of species, the setting in place of no take reserves is often mentioned. Most work on this subject, theoretical as well as applied, was undertaken in recent years. In this paper, we seek to compare two different models presented in existing literature by highlighting their underlying assumptions. Both models were derived from what is often referred to as the "model of Schaefer-Clark'' (reference to the work of the last author on Mathematical Bioeconomics : Clark [7]). We show that various variations of this model lead to properties that can be very different.


%\pagebreak

\section{Introduction}

It is now well known and largely accepted within the scientific community that the exploitation of halieutics resources has reached a critical threshold and solutions must be found urgently to conserve marine biodiversity -- and indeed the existence of certain species.  Beyond the application of quota policies, other strategies have been proposed such as the creation of marine reserves.  The study of the role of reserves in fishery management has been the subject of renewed interest in recent years \cite{arnason,car,conrad,hann,houde,Pezzey}.  The International Conference on the Economics of Marine Protected Areas (MPA) held in July 2000 in Vancouver was one of the starting points in the development of this new paradigm: the use of MPAs as an instrument in the management of fisheries.

In economics literature, Sanchirico and Wilen \cite{san2001} seem to have been the first to suggest that MPAs could be beneficial not only from an ecological but also from an economic point of view.  In their dynamic and spatial model of a Marine Reserve Creation, they analyze whether the transfer of biomass from the reserve to areas where catch is allowed could create economics profit from the MPA since its creation could help to improve a depleted biomass and increase catch outside of the reserve.  They call this a double-payoff because in this case the MPA would increase both biomass and economic profits from the fishery.  
 
Both theoretical and applied aspects of the subject are well documented in the literature. From a theoretical point of view, numerous types of models have been proposed based on differential equations, mixing an inter temporal dynamic that corresponds to the population growth under consideration with a spatial distribution of this dynamic over diverse zones \cite{dubey,san,san2001}.

Sanchirico and Wilen construct a model where "the population structure
is characterized in a manner consistent with modern biological ideas
that stress patchiness, heterogeneity and interconnections among and
between patches"  \cite{san99}.  In this model, independent growth dynamics thus are associated with different patches.

Another type of model has been used in the literature to account for an analogue structure \cite{sri} where a population develops different characteristics in sub-zones.  In these models, the population of the entire zone under consideration follows a given dynamic evolution and the diverse sub-zones have dynamics such that by aggregating them together one may rediscover the global dynamic.

The question that one then must ask is whether these two approaches
may be used interchangeably.  Few studies have focused on comparing
these different model types.  We would like to demonstrate that the choice between a patch model and a global model is not a neutral one, and highlight characteristics of these two models that often are not specified in the literature.

The two models that we will compare are both Clark type (Schaefer, Gordon,...) \cite{clark} whose dynamic is logistic (Verhulst) and therefore widely used in halieutics dynamics.  We use them to study the consequences of setting up a marine reserve from both an economic (inter temporal revenue) and a biologic point of view (population stocks).

The first, the patch model, assumes a relative autonomy between reserve and non-reserve zones.  This type of model is fairly widely used.  The second model, known as the global model, assumes for its part a greater interaction between zones.

The rest of the paper is organized as follows.  In Section 3, we introduce and examine the two models we wish to study.  We then compare results obtained, particularly using numerical simulation, in Section 4.  Section 5 concludes and is followed by a series of annexes that present the demonstration of various results. 

\section{Two variations on the Clark model}
We present two modelizations for a
protected area in a given zone. These two modelizations derive from
the well known fishery model studied for instance by Clark
\cite{clark} among others.
\subsection{ The first variation:  patches model}
\subsubsection{The model}
We consider a fish population that lives in a zone caracterised by a carrying
capacity $K=1$. We assume that this zone splits in two sub-zones, with
capacity respectively equal to $\alpha$ and $ 1-\alpha$; we denote the
stock of the
corresponding sub-populations by $x_1$ and $x_2$.

These two populations follow two independent evolutions laws and this is
the reason why we use the ``patches-concept''. These evolutions
are given by:

\begin{equation}
 {d x_1 \over d t}(t):= \dot x_1(t)=F_{1}(x_1(t))
 \nonumber
\end{equation}
\begin{equation}
 {d x_2 \over d t}(t):= \dot x_2(t)=F_2(x_2(t)).
 \nonumber
\end{equation}

The standard reference for the evolutions law is the logistic law:
$$F_1(x_1):= r_1 x_1(1-{x_1 \over \alpha})$$ 
$$F_2(x_2):= r_2 x_2(1-{x_2 \over 1-\alpha}).$$

But for our purpose it is enough to assume that the $F_i$ are strictly
concave, $C^1$ functions defined on
 $[0,\alpha]$
respectively on $[0,1-\alpha]$ and the $F_i$
 satisfy  $F_i(0)=0, F_1(\alpha)=F_2(1-\alpha)=0$. \\

We now assume that some exchange exists between these two patches and
that this can be represented in terms of the density of these
populations. More precisely we assume that the existence of some diffusion  between these two
patches can be captured by the following:
\begin{equation}
\lambda ({x_2 \over 1-\alpha}-{x_1 \over \alpha})
\nonumber
\end{equation}
where  $\lambda \geq 0 $ represents a diffusion coefficient (${x_2 \over
  1-\alpha}$, ${x_1 \over \alpha} $ being  the density of the
  populations). The value of the diffusion coefficient depends on the
  location of the protected area.

\vspace{1cm}

From now on, we decide that the first zone with capacity $\alpha$ is a protected
area where no catch is allowed, whereas in the second zone
fishing is allowed.\\
The ``normal'' situation, i.e.  the protected
area acts like a source of biomass, corresponds to the case where the
density inside the protected area is bigger than outside, i.e. 
$${x_2 \over 1-\alpha} \leq {x_1 \over \alpha}. $$ 
The growth of the two sub-populations are governed respectively
by the following dynamics:
\begin{eqnarray}
\begin{array}{rrl}
\dot x_1(t)&=&F_1(x_1(t))+\lambda (t)({x_2(t) \over 1-\alpha}-{x_1(t) \over \alpha})\\
\dot x_2(t)&=&F_2(x_2(t))-\lambda (t)({x_2 (t)\over
  1-\alpha}-{x_1(t) \over \alpha})-h(t)
\label{contraintes}
\end{array}
\end{eqnarray}
where  $h(t)$  is the capture rate at time $t$. \\
We note, from the positiveness of $\lambda $ and of the functions
$F_i(.)$, that if the sytem (\ref{contraintes}) possesses an
equilibrium, it has to be necessarily normal.

As it is generally assumed the catch is proportional to the fishing
effort $E$, and to the density of the  population \cite{boncoeur}, therefore
given by:
$$h(t)=  qE(t){x_2(t) \over 1-\alpha} $$
where $q$, the catchability coefficient, represents the fishing death
rate when the density of the population is equal to one. We assume 
 \begin{equation}
0\leq E(t)\leq E_M \quad \quad q>0
 \nonumber
\end{equation}

\vspace{1cm}

The catch is sold on a market. In order to simplify we assume a
constant 
price, $p$, over time and a constant cost, $c$, proportional to the effort.
Therefore the revenue  at $t$ time is

\begin{equation}
ph(t) -cE(t)=(pq{x_2(t) \over 1-\alpha} - c) E(t)
 \nonumber
\end{equation}

We then consider the discounted total revenue on an infinite horizon is 
given by
\begin{eqnarray}
J(E(.),\lambda(.)):= \int_0^{\infty}\; e^{-\delta t}(pq{x_2(t) \over 1-\alpha}-c)E(t) 
\;dt
\end{eqnarray}
where $\delta>0$ is an actualisation factor.

We assume the existence of a manager whose goal is the
maximisation of this total revenue. Moreover we assume that this
manager can act on the
fishing effort $E$ and on some caracteristics of the reserve (closure,
location) that are captured by 
$\lambda $.  
Then the manager
  faces the following control problem:

\begin{eqnarray}
\begin{array}{rl}
\max\limits_{E(.),\;\lambda(.)}& J(E(.),\lambda(.)) \\
\mbox{s.t.}& (\ref{contraintes}) \\
 & 0\leq E(t)\leq E_M
\label{Pb1}
\end{array}
\end{eqnarray}

{\bf Remark } In many papers the states variables stand for the  densities of the
populations and not for the amount of the biomass. The link with the
present model is obtained in setting:
$$X_1={x_1 \over \alpha}, \;\; X_2={x_2 \over 1-\alpha} $$
The two dynamical equations that give the evolution of the populations in the logistic case,  then become:
\begin{eqnarray}
\begin{array}{lll}
 \dot X_1& =& r_1X_1(1-X_1) + {\lambda \over \alpha}(X_2-X_1) \\
 \dot X_2& =& r_2X_2(1-X_2) - {\lambda \over 1-\alpha}(X_2-X_1)- Q E X_2
\end{array}
\end{eqnarray}
with $Q={q \over 1-\alpha} $.\\ This model with patches could be considered as a more
general one than the  \cite{Ami} paper which corresponds to  $\alpha = {1 \over 2}$.

\subsubsection{Analysis of the solutions}
We will study the previous optimal control problem by
the help of the calculus of variations theory .\\
From the dynamic (\ref{contraintes}), we deduce the expression of the
effort in terms of the state  variables:

%\begin{equation}
$$ E(t) =\frac{1-\alpha }{q x_2(t)}(F_1(x_1(t))+F_2(x_2(t))-\dot 
x_1(t)-\dot x_2(t) )$$
% \label{effort}
%\end{equation}
and then we obtain the new form of the objective. Thus the optimisation problem  becomes:
%\begin{eqnarray}
%\begin{array}{l}
$$\max\limits_{X\in C}\int_0^{\infty}\; e^{-\delta t}(p-{c (1-\alpha) \over
  qx_2})[F_1(x_1)+F_2(x_2)-\dot x_1-\dot x_2] \;dt$$
%\\
%\end{array}
%%\label{calvar}
%\end{eqnarray}

$C$ being the set of admissible curves:
\begin{eqnarray}
\begin{array}{rl}
C= \{x(.)=(x_1(.),x_2(.)) \;\;& x_i(.) \in BC^1([0,\infty[),x_i(0) \mbox{ given },\\
 &G(x_1,x_2)-qE_M{x_2 \over 1-\alpha} \leq \dot x_1+\dot x_2 \leq G(x_1,x_2)\}  
\end{array}
\end{eqnarray}
with  \; $ G(x_1,x_2):= F_1(x_1)+F_2(x_2)$ and $BC^1$ stands for the
bounded with bounded derivative functions defined on the interval $[0,\infty[$.
 
It is known that on $BC^1$,
the first order optimality conditions given by the Euler-Lagrange
equations, apply (see \cite{Blot}). We suppose that $x(.)$ stands for an interior solution and then $x(.)$ has to satisfy:

 $$l_{x_i}(x_1(t),x_2(t))-\frac{d}{dt}l_{\dot 
x_i}(x_1(t),x_2(t))+
 \delta l_{\dot x_i}(x_1(t),x_2(t))=0$$
$l(.,.)$ being the non actualised Lagrangian of the calculus of
variations problem, $l_{x_i}(.,.)$ stands for the derivative with
respect to $x_i$ and   $l_{x_i}(.,.)$ stands for the derivative with
respect to $x_i$..
\\
The Euler-Lagrange equations becomes then

\begin{eqnarray}
\begin{array}{rrl}
\dot x_1&=&x_2({pqx_2 \over c(1-\alpha)}-1)(F_2'(x_2)-\delta)+
F_1(x_1)+ F_2(x_2)\\
\dot x_2&=&x_2({pqx_2 \over c(1-\alpha)}-1)(\delta-F_1'(x_1)).
\end{array}
\label{Euler}
\end{eqnarray}
We first are interested by the non trivial equilibria ,
$(x_1^*,x_2^*)$,  of (\ref{Euler})  i.e. such that $F_i(x_i^*) \neq
0, $ i.e. $x_1^* \neq 0,  \alpha$ and $x_2^* \neq 0, 1-\alpha $.
It is easy to establish that such equilibria have to satisfy (Appendix 1)
\begin{eqnarray}
 F_1'(x_1^*)=\delta.
\nonumber
\end{eqnarray}

We now assume that $r_1> \delta $. Then in the logistic case, from the strict concavity
 of  $F_1(.)$, we immediatly deduce the existence of a unique $x_1^* \in ]0,\alpha/2[$. \\ Therefore we obtain the
following result whose proof is postponed in Appendix 1.

\vspace{1cm}

\noindent
{\it {\bf Lemma 1} In the logistic case with $r_1 > \delta $,  there is a unique positive
non trivial solution, $(x_1^*,x_2^*)$, of the Euler-Lagrange equations
(\ref{Euler}). This solution is caracterised by  $$x_1^*= {\alpha (r_1-\delta)
  \over 2r_1}$$ and $x_2^*$ given by 
\begin{eqnarray}
x_2[{2r_2 pq \over c(1-\alpha)^2}x_2^2-({pq \over c(1-\alpha)}(r_2-\delta)+{r_2 \over
  1-\alpha})x_2 -\delta
]=\alpha {(r_1-\delta)(r_1+\delta) \over 4r_1}.
 \label{etat2}
\end{eqnarray} }

\vspace{1cm}

Clearly  $x_1^* \in ]0,\alpha[$. It remains to show  that  $x_2^*
\in [0,1-\alpha]$. This can be done straightforwordly but we prefer to
use the following approach. We recall that if (\ref{contraintes}) possesses
a non-zero equilibrium $(x_1^*,x_2^*)$ then it is necessarily normal that is to
say 
 ${x_2^* \over 1- \alpha} \leq {x_1^* \over \alpha}$ and we   first observe
 that if this last condition holds then $x_2^* \leq {x_1^* \over
   \alpha}(1-\alpha) \leq 1-\alpha$. \\
We prove now that this condition holds under conditions on
the ratio ${pq \over c}$ and on the size of the different zones given
by $\alpha $.
\\ Let us introduce 
 $T(.)$ defined from the left hand member in (\ref{etat2})
%\begin{eqnarray}
$$ T(z) = z [2r_2 \theta z^2-(\theta (r_2-\delta)+r_2)z -\delta ]$$
%\label{theta}
%\end{eqnarray}
where $\theta:={pq \over c}$.
Then:

$$T({x_2^* \over 1- \alpha})= {\alpha \over 1-\alpha} {(r_1-\delta)(r_1+\delta) \over 4r_1}.$$

As we assume $r_1> \delta $, then $T({x_2^* \over 1-\alpha}) >0$.
 From the graph of $T(.)$ we deduce that the inequality $${x_2^* \over 1- \alpha} \leq
{x_1^* \over \alpha}$$ holds if and only if 
 $$T({x_2^* \over 1-\alpha})\leq
T({x_1^* \over \alpha}).$$ 

Let us introduce when $r_1\neq r_2$ 
\begin{eqnarray}
 \theta_0:=({2r_1 \over r_1-\delta} + {r_2 \over \delta}){r_1 \over
     r_1-r_2}.
\label{teta0}
\end{eqnarray}

The following result is detailled in
Appendix 2.

\vspace{1cm}

\noindent
{\it 
{\bf Lemma 2}\\ In the logistic case,
\begin{enumerate}
\item if $r_1>r_2$,
for each $ \theta > \theta_0 $ if $\alpha$ satisfies
\begin{eqnarray}
\alpha (r_1+\delta + {r_1-\delta \over r_1}(\theta \delta {r_1-r_2
  \over r_1}-r_2)-2\delta) \leq 
 {r_1-\delta \over r_1}(\theta \delta {r_1-r_2 \over r_1}-r_2)-2\delta \ \ \ 
\label{alpha0}
\end{eqnarray}
then the solution $(x_1^*,x_2^*)$, of the Euler Lagrange equation
(\ref{Euler}) given in Lemma 1, is normal i.e. satisfies:
$${x_2^* \over 1- \alpha} \leq {x_1^* \over \alpha}$$
\item 
if $r_1\leq r_2$,  $(x_1^*,x_2^*)$ is never normal i.e.
$${x_2^* \over 1- \alpha} > {x_1^* \over \alpha}.$$
\end{enumerate}
}

\vspace{1cm}

Then it remains to prove that $(x_1^*,x_2^*)$ is a candidate  to the optimisation
problem (\ref{Pb1}), thus that the constraints are satisfied.  

If the conditions in the Lemma 1 and Lemma 2.1 hold, then  $F_1(x_1^*) >
0$ and  $\displaystyle {x_2^* \over 1-\alpha} < {x_1^* \over
  \alpha^*}$. Therefore to this corresponds a unique coefficient of
diffusion  $\lambda^* > 0$.
\\ From (\ref{contraintes}) we deduce that
\begin{eqnarray}
E^*={1-\alpha \over q x_2^*}(F_1(x_1^*)+F_2(x_2^*)) > 0.
\label{effort1}
\end{eqnarray}
  Moreover
the expression of the total revenue is given by
\begin{eqnarray}
 J^* = 
(pq {x_2^* \over 1-\alpha}-c){ E^* \over \delta}.   
\label{Revenu1}
\end{eqnarray}
This revenue is positive if the fishery profit is positive at this
equilibrium, that is to say if 
$$\displaystyle {x_2^* \over 1-\alpha} > {c \over pq}={1 \over \theta}.$$
But this inequality holds because we have 
$$T({x_2^* \over 1-\alpha})>T({1 \over \theta})= (1-\theta){r_2
  \over \theta ^2}$$
and due to the fact that $\theta \geq  \theta_0 > 1$,
 this last term is nonpositive .

\vspace{1cm}

\noindent
{\it {\bf Proposition 1} \\
In the logistic case if $r_1>\delta$ the problem (\ref{Pb1}) possesses at most a non trivial and positive
optimal stationary solution caracterised  by 
$$x_1^*= {\alpha (r_1-\delta) \over 2r_1}$$ and $x_{2}^{*}$ given by (\ref{etat2}).
The corresponding effort, diffusion coefficient and total revenue are
given respectively by (\ref{effort1}), (\ref{contraintes}), (\ref{Revenu1}).
}

{\bf Remarks}\\
1) When $\alpha ={1 \over 2}$ the value of  $\theta _0 $ coincides with
the value $\tilde p_m$ given in \cite{Ami}.\\
2) From the expression of $\theta _0$ in (\ref{teta0}), we observe that
$r_1$ can't be closed to $r_2$. If this is not the case, then the value of the dimensionless
ratio $\theta$ has to be very high. But this can be unrealistic
because the value of $\theta$ is given by the economic environment. \\
3) From (\ref{teta0} ) with a given value for $\theta$ we can precise
a bound for $r_2$ expressed in terms of $r_1, \delta$:
$$r_2\leq  {\theta -{2r_1\over r_1 -\delta} \over \theta + {r_1\over
    \delta}} r_1.    $$

\subsection{The second variation: the splitting of a unique zone}
In this second model we start with an unique zone with capacity $K$ that
we normalise to one, $K=1$. Let us assume that the population follows
a standard  evolution law:
\begin{eqnarray} 
\dot z(t)=\phi (z(t)) \quad ( \mbox{for instance}\quad  = rz(t)(1-z))
\nonumber
\end{eqnarray}
 $\phi (.)$  being a $C^1$ concave function defined on
$[0,1]$,  with
$\phi(0)=\phi(1)=0$ and $\phi'(1) <0$.
We assume that this zone splits first in a part that is a reserved area where
no fishing is allowed and a complementary part that is open to harvest. We
assume that these two parts have respectively $\alpha$ and  $1-\alpha$
as a carying capacity. \\
The main difference with the previous model is that the two
populations, which stocks are respectively  $x_1$ and $x_2$ follow the 
evolution laws given by:
\begin{eqnarray} 
\dot x_1(t)=F(x_1(t),x_1(t)+x_2(t)) \nonumber \\
\dot x_2(t)=F(x_2(t),x_1(t)+x_2(t)) \nonumber
\end{eqnarray}
where $F(.,.)$ satisfies the standard assumption of regularity with $ F(x,z)=0 \quad \mbox{ if and only if } \quad x= 0 
\mbox{ or  } 
z=1$
and where $F(.)$ and  $\phi(.)$ satisfy 
\begin{eqnarray}
F(x_1,x_1+x_2)+F(x_2,x_1+x_2)=\phi (x_1+x_2). \nonumber
\end{eqnarray}
For instance $F$ can be a logistic function $$ F(x,z)=
rx(1-z)=rx(1-(x_1+x_2)).$$
As in the previous model, there is  some diffusion between the two
zones which can be
represented by:
$$\lambda ({x_2 \over 1-\alpha}-{x_1 \over \alpha}).$$
Then the two populations evolve following the dynamics:
\begin{eqnarray} 
\dot x_1(t)=F(x_1(t),x_1(t)+x_2(t))+ \lambda ({x_2(t) \over 1-\alpha}-{x_1(t) \over 
\alpha}) \nonumber\\
\dot x_2(t)=F(x_2(t),x_1(t)+x_2(t))-\lambda ({x_2(t) \over 1-\alpha}-{x_1(t) \over
  \alpha}).
\nonumber
\end{eqnarray}
We want to stress on the fact that this new model is  consistant in the sense that  the
sum of the two dynamics is exactly the evolution law of the total population.

Now taking into account the catch in the zone where fishing is allowed,
we derive the final dynamics of the populations

\begin{eqnarray}
\begin{array}{rrl} 
\dot x_1(t)&=&F(x_1(t),x_1(t)+x_2(t))+ \lambda ({x_2(t) \over 1-\alpha}-{x_1(t) \over 
\alpha})\\
\dot x_2(t)&=&F(x_2(t),x_1(t)+x_2(t))-\lambda ({x_2(t) \over 1-\alpha}-{x_1(t) \over
  \alpha})-qE{x_2(t) \over 1- \alpha }
\end{array}
\label{dyn2}
\end{eqnarray}
 where $E\in [0,E_M] $ stands for the fishing effort and $q>0$ is the catchability coefficient.

In order to compare with the previous model, we assume that a manager
has as an  objective to maximise  the actualised total revenue
as presented before. To do so, he has to act on two controls, the fishing effort and the
location of the reserve area given by $\lambda$. Therefore the manager faces the following
program of optimisation
\begin{eqnarray}
\begin{array}{rl}
\max\limits_{E(.),\; \lambda(.)}&\int_0^{\infty}\; e^{-\delta
  t}(pq{x_2(t) \over
 1- \alpha }-c)E(t) \;dt \\
 \mbox{s.c.} &(\ref{dyn2}) \\
&0\leq E(t) \leq E_M \\
\label{pb2}
\end{array}
\end{eqnarray}

\subsubsection{Analysis of the solutions}
From the dynamic equations (\ref{dyn2})  we can derive the expression
of the effort
\begin{eqnarray}
E(t)= {1-\alpha \over qx_2(t)}(\phi (x_1(t)+x_2(t))-\dot x_1(t)-\dot x_2(t))
\label{effort2}
\end{eqnarray}
 and then we obtain the equivalent problem to
(\ref{pb2}) as a 
calculus of variations problem:
%\begin{eqnarray}
$$\max\limits_{X\in C}\int_0^{\infty}\; e^{-\delta t}(p-{c(1-\alpha) \over qx_2})[\phi
(x_1+x_2)-\dot x_1-\dot x_2] \;dt$$
%\label{calvar2}
%\end{eqnarray}

where  $C$ stands for the set of feasible curves defined by:
\begin{eqnarray}
\begin{array}{rl}
$$C= \{x(.)=(x_1(.),x_2(.)) \;\;& x_i(.) \in BC^1([0,\infty[),x_i(0) \mbox{ given },\\
& \phi(x_1+x_2)-qE_M{x_2 \over 1-\alpha} \leq \dot x_1+\dot x_2 \leq \phi(x_1+x_2)\}. 
%\label{C2}
 \nonumber
\end{array}
\end{eqnarray}

We know that in this framework a necessary optimality condition for an interior solution $x(.)$ is given by the Euler-Lagrange
equations that are:
\begin{eqnarray}
\begin{array}{rrl}
\dot x_1&=&x_2({pqx_2 \over c(1-\alpha)}-1)(\phi '(x_1+x_2)-\delta)+
\phi(x_1+x_2)\\
\dot x_2&=&x_2({pqx_2 \over c(1-\alpha)}-1)(\delta-\phi '(x_1+x_2)).
\end{array}
%\label{Euler2}
 \nonumber
\end{eqnarray}
For now, we will stick to the logistic case. The Euler-Lagrange
equations are then:
\begin{eqnarray}
\begin{array}{rrl}
\dot x_1&=&x_2({pqx_2 \over c(1-\alpha)}-1)(r-2r(x_1+x_2)-\delta)+
r(x_1+x_2)(1-(x_1+x_2))\\
\dot x_2&=&x_2({pqx_2 \over c(1-\alpha)}-1)(\delta-r+2r(x_1+x_2)).
\end{array}
\label{Euler3}
\end{eqnarray}

\vspace{1cm}

In order to derive the non trivial positive equilibria, denoted by $(x_1^*,x_2^*)$, we
first consider the second equation in (\ref{Euler3}) with  the condition  $$x_1^*+x_2^*={r-\delta \over 2r}.$$ 
This implies as a result
$$0=
r(x_1^*+x_2^*)(1-(x_1^*+x_2^*))={(r-\delta)(r+\delta)
  \over 4r}.$$ 
A contradiction if $r>\delta$. In the case where $r=\delta$, then we
obtain the trivial solution $x_1^*=x_2^*=0$. \\ Then we deduce that an equilibrium
has  to necessarily satisfy
$$x_2^* = {c(1-\alpha) \over pq}.$$
With the help of first equation in (\ref{Euler3}) we find that either
$x_1^*+x_2^*=0$ or  $x_1^*+x_2^*=1$. \\
Finding a non trivial equilibrium implies to exclude the first condition. 
Therefore we have proved the following result

\vspace{1cm}

\noindent
{\it 
{\bf Lemma 3} In the logistic case there is a unique non trivial and positive solution
for the Euler-Lagrange solutions (\ref{Euler3}) given by: 
$$(x_1^*,x_2^*)=(1-{c(1-\alpha) \over pq},{c(1-\alpha) \over pq}).$$}

\vspace{1cm}

In order to examine whether  this candidate solution of the
problem (\ref{pb2}) can be optimal or not, we have to derive the corresponding effort and
diffusion coefficient. From the expression of the effort
(\ref{effort2}), we obtain   $\phi(x_1^*+x_2^*)=\phi(1)=0$
that is to say
$$E^* =0.$$
We also deduce that
$$\lambda ^* =0$$
 except if  $\theta =1$.
Finally at this equilibrium  the intertemporal revenue is null too.
We have established the following proposition
\vspace{1cm}

\noindent
{\it 
{\bf Proposition 2} \\
In the logistic case the problem (\ref{pb2}) possesses at most a non trivial and positive
stationary solution given by
$$(x_1^*,x_2^*)=(1-{c(1-\alpha) \over pq},{c (1-\alpha)\over pq}).$$
The corresponding effort, diffusion coefficient and total revenue are
null.
}

\noindent
{\bf Remarks}\\

1) It is easy to obtain that this equilibrium is normal, i.e.
$${x_1^* \over \alpha}  > {x_2^* \over 1-\alpha}$$
if the fishery is profitable, that is to say if $$pq-c>0.$$
\\

2) An adaptation of  the model given in Gomez et al. \cite{cartigny} to our case of a no take zone
is: 
\begin{eqnarray}
\begin{array}{rrl}
\dot x_1 & = & \alpha r(x_1+x_2)(1-(x_1+x_2)) +\lambda ({x_2 \over
  1-\alpha}-{x_1 \over \alpha}) \\
\dot x_2 & = & (1-\alpha)r(x_1+x_2)(1-(x_1+x_2)) -\lambda ({x_2 \over
  1-\alpha}-{x_1 \over \alpha}) - q_{2}E_{2}{x_2 \over
  1-\alpha}.
\end{array}
\end{eqnarray}
 We can get the same results as those given earlier. In Gomez et
 al. \cite{cartigny},  fishing is allowed in the so called artisanal
 zone (corresponding to the protected area in our case) and the objective to be maximised is somewhat
different (it takes into account the revenues of the artisanal and
industrial fisheries). Here also it has been proved that a unique solution  exists
 but with a non null effort and a non null revenue.

\section {Comparison, Numerical application} 
\noindent
In this section we underline the differences between the results we
obtained in the previous sections for both the patches case and the
global model. \\  
From their expressions given in the Propositions 1 and Proposition 2,
we can make the following remarks for the equilibria $(x_1^*,x_2^*)$: 

\begin{itemize}
\item 
In the model with patches the first component $x_1^*$ doesn't depend explicitely
on the ratio ${c\over pq}$, whereas it does in the global model.
\item In the global model the second component is given by $x_2^*={c(1-\alpha) \over
  pq}$, whilst the patches model doesn't possess any equilibrium with
such a component as the second, cf. Appendix 1.
\end{itemize}
\noindent
Thus the expressions of the equilibria are different in the two models. 

\vspace{0.5cm}
 
Now we established that the optimal effort and the corresponding total
revenu at $(x_1^*,x_2^*)$ was null for the global model. This doesn't
seem to be the case for the patches model, we will later show with
simulations that optimal effort and total revenue are not
significantly close to zero. 

\vspace{0.5cm}
 
In order to continue the comparison, let's arbitrarily fix
the parameters
$p, q, c, \alpha, \delta $. Thus the models depend only on the instantaneous
growth rates $r_1, r_2$ and $r$ respectively. \\
If we let $r_1=r_2=r$, in Lemma 2 we established that the equilibrium
$(x_1^*,x_2^*)$ was never normal in the patches case, while it is
always normal for the global model (Remark 1 in section 3.2). 
\\
Now to compare our models with $r_1$ and $r_2$ only near $r$, we noted
in the Remark 2 of section 3.1 that this situation wasn't a realistic one.

Then comparing these two models is not  an easy task. Our first
conclusion is:   the role
of the instantaneous growth rates of the biomasses are crucial to
choose such or such model.  An
assumption that is not underlined in general.

Now let's come back to the comparison of the optimal efforts and revenues by
using simulations. The main issue is to determine significant
growth rates that are not equal.
\\
But this choice shouldn't depend on our particular models with 
preserving areas. It should be the same for a wide class of models.  For
instance for models that correspond to a situation where fishing is allowed in the
the two areas (\cite {Raissi}). We consider thus:
\begin{eqnarray}
\begin{array}{lll}
\dot x_1& =& F_1(x_1) + \lambda({x_2 \over 1- \alpha} -{x_1 \over \alpha}) - qE x_1 \\

\dot x_2 &=& F_2(x_2) - \lambda({x_2 \over 1- \alpha} -{x_1 \over \alpha}) - qE x_2
\label{1'}
\end{array}
\end{eqnarray}
and
\begin{eqnarray}
\begin{array}{lll} 
\dot x_1& =& F(x_1,x_1+x_2) + \lambda({x_2 \over 1- \alpha} -{x_1 \over \alpha}) - qE x_1\\

\dot x_2 &=& F(x_2,x_1+x_2) - \lambda({x_2 \over 1- \alpha} -{x_1 \over \alpha}) - qE x_2
\label{2'}
\end{array}
\end{eqnarray}
with the same assumptions as before.
In order to compare numerically (\ref{1'}) and  (\ref{2'}) we will 
face the same issue to determine significant growth rates. 
\\ We propose to use this new situation in order to fix values for
$r_1,r_2,r$. The new problem we consider now is to maximize the same
objective as before  
 $$J(E) = \int_{0}^{\infty} (pq(x_1+x_2)-c) Ee^{-\delta t} dt$$ 
but with  (\ref{1'}) and  (\ref{2'}). 
\noindent
We observe that  (\ref{2'}) corresponds to the classic Clark
model, it is enough to let
$z=x_1+x_2$ to obtain that the dynamic is $\dot z =\phi (z) -qEz$. 

\vspace{0.5cm}

We can assume that a manager has no reason to use one model rather
than another. Then the two models can be considered as equivalent in the sense that
they provide the same optimal effort.\\ Hence, we propose the following
procedure to determine a system of growth rates: 
Let's set an arbitrary choice of values for $r_1$ and $r_2$. From the
first order
optimality conditions, given here by the Pontryagin principle, we can
derive the optimal value of the corresponding effort for problem
(\ref{1'}). We hand-over this value in the first order optimality
conditions of the second problem  (\ref{2'}) from which  we derive the value
of the growth rate $r$.
\\ To follow this procedure we set:
$$ \alpha=.5, \delta=.05, c=.15, q= 2., \lambda=20$$
and we obtained 
for
$$r_1=0.4, r_2 =0.05 $$ that  $$\bar E= .0566 \;\; \mbox{ and } \; \; r=0.28739 .$$

\vspace{0.5cm}

Let's  now go back to our models with protected areas from where we take the
previous values for the parameters and where we set $r_1=.4, r_2 =.05 $ and $ r=0.28739$ . 

\noindent
Then for the model with patches we found that:
\begin{itemize}
\item 
the optimal effort is  $E^*=.0457$, and the biomass values are
respectively   $x_1^*=.21875, x_2^* = .0302$
\end{itemize}
and for the global model
\begin{itemize}
\item 
 $E^*=0$ et $x_1^*=.875, x_2^* = .125$ .
\end{itemize}

We observe that the optimal efforts corresponding to the patches case,
$E^*=.0457$, and the Clark model (\ref{2'}), $\bar E=0.0566 $,
have similar sizes. We know that the optimal value of the effort in
this last model can't be
considered as null. Therefore we can  deduce that in the first model
with patches the effort is not null. \\ Then  the two models have
different qualitative behaviour:total revenues and  optimal effort are
totally different.

\section{Conclusion}

In this article, we have shown that different models have been proposed and used in the literature for the same MPA problematic. We focused on demonstrating the importance of the hypotheses underlying two types of models -- the patch model and sub-zone model --particularly the crucial role played by the growth functions (rate and form), and on studying the different results produced by them.

A manager who wishes to study the role of an MPA in a given zone must first know if the entire zone is artificially divided or if it can be broken down into patches (entities with their own dynamics).  Without taking this precaution, and in obtaining the very different results that we have seen, the manager risks taking erroneous decisions.

The two preceding models of resource dynamics are adapted to the case where control instruments are independent of the size of the no take reserve. If the manager must take size into account in his decisions, the modelling of the
 dynamic has to be changed. For instance, it is necessary to consider
 a depending of $\alpha $ diffusion coefficient. A justification is
 given in Appendix 3. In this Appendix, we also underline that this
 coefficient can be given by formula
 $$ \lambda (\alpha)= \lambda _0 \alpha (1-\alpha)$$
 which is the expression considered by Boncoeur (cf. \cite{boncoeur}).

\section{Appendix}
\subsection{Appendix 1}
We determine  the non trivial equilibria of the Euler-Lagrange equations
(\ref{Euler}) 
\begin{eqnarray}
\begin{array}{rrl}
\dot x_1&=&x_2({pqx_2 \over c(1-\alpha)}-1)(F_2'(x_2)-\delta)+
F_1(x_1)+ F_2(x_2)\\
\dot x_2&=&x_2({pqx_2 \over c(1-\alpha)}-1)(\delta-F_1'(x_1)).
\end{array}
\nonumber 
\end{eqnarray}
If we assume that $x_2^*={c(1-\alpha) \over pq}$, from the first
equation:
\begin{eqnarray} 
F_1(x_1)+F_2(x_2^*)=0.
\label{Sum}
\end{eqnarray}
From the assumption of the non triviality of the equilibria, we have that $c\neq 0,
\alpha \neq 1, c\neq pq$ and then $F_2(x_2^*)>0$. Therefore we can't
find 
any $x_1\in [0,\alpha]$ such that (\ref{Sum}) holds. Thus there is no non
trivial equilibrium with $x_2^*={c(1-\alpha) \over pq}$.

Therefore in order for a non trivial equilibrium to exist it is necessary that 

$$ F'_1(x_1^*))=\delta. $$

In the logistic case it is easy to compute that 
 $$x_1^*= {\alpha (r_1-\delta) \over 2r_1}.$$ 
Reporting this value in (\ref{Euler})
we get  the following equation for $x_2^*$
\begin{eqnarray}
x_2({pqx_2 \over c(1-\alpha)}-1)(F_2'(x_2)-\delta)+
F_2(x_2)=-F_1(x_1^*)
\nonumber
\end{eqnarray}
which is (\ref{etat2}) in the logistic case
\begin{eqnarray}
x_2[{2r_2 pq \over c(1-\alpha)^2}x_2^2-({pq \over c(1-\alpha)}(r_2-\delta)+{r_2 \over
  1-\alpha})x_2 -\delta
]=\alpha {(r_1-\delta)(r_1+\delta) \over 4r_1}. \nonumber
\end{eqnarray}
The graphes of the functions defined by the left and right hand members are curves that crosse in
a single $x_2^*$  if $r_1\geq \delta$. But this last
inequality holds because from our assumption we have
\begin{eqnarray}  
 F_1(x_1^*)=\alpha {(r_1-\delta)(r_1+\delta) \over 4r_1} \geq 0.
\nonumber
\end{eqnarray}
This ends the proof of Lemma 1.

\subsection{Appendix 2}

In order to find conditions for the inequality ${x_2^*\over 1- \alpha}
\leq {x_1^* \over \alpha}$ to be true, we know that it is equivalent
to consider the
inequality  $T({x_2^*\over 1- \alpha})
\leq T({x_1^* \over \alpha})$.
This last inequality becomes
\begin{eqnarray}
\begin{array}{l}
{\alpha \over 1-\alpha} {(r_1-\delta)(r_1+\delta) \over 4r_1} \leq \nonumber \\T({x_1^* \over \alpha})=
{r_1-\delta \over 2 r_1} [{2r_2 pq \over c}({r_1-\delta \over 2 r_1})^2
-({pq \over c}(r_2-\delta)+r_2){r_1-\delta \over 2 r_1} -\delta ]
\nonumber
\end{array}
\end{eqnarray}
that is equivalent to
\begin{eqnarray}
 \alpha (r_1+\delta)\leq (1-\alpha)({r_1-\delta \over
  r_1}(\theta \delta {r_1-r_2 \over r_1}-r_2) -2\delta). 
\nonumber  
\end{eqnarray}

1) If $r_1> r_2$
the right hand member has to be positive, this implies the
following condition on $\theta$:
\begin{eqnarray}
 \theta > \theta_0:=({2r_1 \over r_1-\delta} + {r_2 \over \delta}){r_1 \over
     r_1-r_2}. 
\nonumber
\end{eqnarray}
Now if this condition on $\theta $
holds, from the previous inequality we should deduce (\ref{alpha0}).

2) If $r_1\leq r_2$, the right hand  member is always negative and
therefore $(x_1^*,x_2^*)$ can't be normal.

This ends the proof of Lemma 2.

\subsection{Appendix 3 }
We consider the case where the manager has the size of the
preserved area as control. We will first prove that the diffusion coefficient
has to depend on this size.
\\
We  start with the dynamics and the objective given in
the second variation  (§ 3.2). We suppose that the manager  has to maximise his objective by using 
the fishing effort $E$ and the size of the preserved area that is
captured by $\alpha$.
\\
We always denote by $z$  the stock of the total population and  the two subpopulations stocks by $x_{1}$ and $x_{2}$
respectively. Then  the densities in the two regions are $d_1={x_1 \over \alpha}$ and
 $d_2={x_2 \over 1-\alpha}$. 
\\
When  $\alpha =0$ we can only find a single zone and thus $x_1=0$ and
$x_2=z$. In this case, it is natural to set for the densities: $d_1=0$ and 
$d_2=z$ respectively. Now if $\alpha =1$, it is natural to set: $d_{1}=z$  $d_2=0$.
\\ As we have done before, we assume that some diffusion exists between the
two zones and that it is proportional to the difference of the two
densities. Therefore in order to respect our previous remark, we have
to set 
\begin{eqnarray}
\lambda (\alpha ) ({x_2 \over 1-\alpha }-{x_1\over \alpha })
\end{eqnarray}
where the diffusion coefficient depends on $\alpha$. Indeed if
$\alpha =0$, from $d_1=0$ we deduce that $\lambda (0) {x_2 \over
  1- \alpha }= \lambda (0)z=0$ because in this case we can only find a
single zone, and
thus $\lambda (0) =0$. From a similar argument, we deduce that for 
$\alpha =1$, we have  $\lambda (1)=0$. Now for $\alpha \in ]0,1[$
the coefficient $\lambda (\alpha)$ is certainly not null. 
 \\ For instance we can let
\begin{eqnarray}
\lambda (\alpha ) =\lambda _0\alpha (1-\alpha )
\label{size}
\end{eqnarray}
and in this case the diffusion is modelised by
\begin{eqnarray}
\lambda _0 \alpha (1-\alpha )({x_2 \over 1-\alpha }-{x_1\over \alpha }).
\end{eqnarray}
This expression is the one proposed by Boncoeur in \cite{boncoeur}.
\\ Then the problem of the manager becomes in this setting
\begin{eqnarray}
\begin{array}{rl}
\max\limits_{E(.),\;\alpha} & \int_0^{\infty}\; e^{-\delta
  t}(pq{x_2(t)\over 1-\alpha}-c)E(t) \;dt \\
	\mbox{s.c.} & (\ref{dyn2})
\end{array}
\end{eqnarray}
where the diffusion coefficient in (\ref{dyn2}) is given by (\ref{size}).
\\ The solution of this problem is straightforward using the
Pontryagin maximum
principle. We won't mention it in this paper.

\end{document}